# Global Solution Branches of Forced Nonlinearly Elastic Single-Director Surfaces

Timothy J. Healey[1]

**Abstract.** We prove the existence of classical solutions for nonlinearly elastic, single-director surfaces or shells characterized by finite mid-plane strains and curvature. We employ topological methods, yielding solution continua far (in norm) from an unloaded, stress-free reference configuration, i.e., a global implicit function theorem of sorts. We first establish the existence of a local solution path, which is subsequently shown to be part of a global solution branch, viz., a maximal connected, locally compact set, emanating from the unloaded state.

## 1. Introduction

We prove the existence of classical solutions for nonlinearly elastic single-director surfaces or shells characterized by finite midplane strains and curvature. The work is motivated, in part, by the problem of wrinkling in highly stretched elastomer sheets, where the incorporation of finite mid-plane strains in the model is crucial, cf. [19], [21]. Of course, direct formulations for such modes are well known, e.g., [2], [20]. Moreover, they are available for use in commercial finite-element codes, e.g., as employed in [21]. The existence of energy minimizers for a class of physically appropriate models of such, allowing for wrinkling, was recently established in [13].

Here we employ topological methods, yielding solution continua "far" (in norm) from an unloaded, stress-free reference configuration – a global implicit function theorem of sorts, cf. Theorems 4.1 and 4.2. The results complement those of [13] and provide a mathematical underpinning for the use of finite-strain shell models in numerical-continuation approaches, e.g., as mentioned above. Also, our results align well with those obtained in the setting of 3D nonlinear elasticity, cf. [10], [14], [15]. However, these do not follow as special cases of those works. While the overall topological approach is the same, the explicit presence of the director field in the stored-energy density here raises new issues. For simplicity, the bounded, stress-free reference configuration is chosen to be flat[2]. The generalization to an initially curved reference configuration is straightforward, especially in view of the coordinate-free notation introduced in [9].

The outline of the work is as follows. We present the formulation of the problem in Section 2. The domain is two-dimensional, while the deformation and the director field, each taking values in $\mathbb{R}^3$, comprise the unknowns. Strong ellipticity is assumed. However, we require a stronger condition at the stress-free reference configuration, accounting for the director as well as the deformation gradient and director-field gradient, cf. (2.12)-(2.16). This represents a generalization of the well-accepted assumption that the elasticity tensor in classical linear elasticity is positive- definite on non-zero symmetric tensors. We assume a parameter-dependent family of "live" body forces and body couples. In Section 3, we prove the existence of a unique, local solution path emanating from the unloaded state via the implicit function theorem. This result already appears to be new. We then employ the nonlinear Fredholm degree of [6] to obtain global continuation results in Section 4. The technical requirements of the former, viz., a Fredholm property and properness of the nonlinear map, are established. The global results of Theorem 4.1 are sharpened in Theorem 4.2 for a special class of stored energies. The argument is based on an idea from [14], employed in the context of 3D nonlinear elasticity.

[1] Department of Mathematics, Cornell University, Ithaca, NY, USA
tjh10@cornell.edu

[2] We prefer the terminology *shell* over *plate*; the difference is irrelevant for finite deformations.

**2. Formulation**

Consider a bounded domain $\Omega \subset \mathbb{R}^2$ with boundary $\partial\Omega$ of class $C^3$, i.e., the boundary is everywhere locally a curve corresponding to the graph of a three-times continuously differentiable function. Throughout, we associate $\mathbb{R}^2 \cong \mathbb{R}^2 \times \{0\} \subset \mathbb{R}^3$, i.e., $\mathbb{R}^2 \cong span\{\mathbf{e}_1, \mathbf{e}_2\}$, with $\{\mathbf{e}_1, \mathbf{e}_2, \mathbf{e}_3\}$ denoting the usual right-handed orthonormal basis for $\mathbb{R}^3$. We associate $\bar{\Omega}$ with the reference configuration of a flat shell or plate. Let $\mathbf{f} : \bar{\Omega} \to \mathbb{R}^3$ denote the deformation and $\mathbf{d} : \bar{\Omega} \to \mathbb{R}^3$ the director field, both presumed sufficiently smooth. Together they determine the configuration of the shell. In the reference state, we have $\mathbf{f} \equiv \mathbf{x} = (x_1, x_2, 0)$ and $\mathbf{d} \equiv \mathbf{e}_3$. Their respective gradients are denoted

$$\nabla\mathbf{f} = \frac{\partial f_j}{\partial x_\mu}\mathbf{e}_j \otimes \mathbf{e}_\mu \text{ and } \nabla\mathbf{d} = \frac{\partial d_j}{\partial x_\mu}\mathbf{e}_j \otimes \mathbf{e}_\mu.$$

Here and throughout, repeated Latin indices sum from 1 to 3, while repeated Greek indices sum from 1 to 2. We assume the existence of a sufficiently smooth stored energy density $W : \Theta \to [0, \infty)$, denoted $W(\mathbf{F}, \mathbf{D}, \mathbf{d})$, where

$$\Theta := \{(\mathbf{F}, \mathbf{D}, \mathbf{d}) : J := \det(\mathbf{F} + \mathbf{d} \otimes \mathbf{e}_3) > \mathbf{0}\}, \tag{2.1}$$

for all $\mathbf{F}, \mathbf{D} \in L(\mathbb{R}^2, \mathbb{R}^3)$ and $\mathbf{d} \in \mathbb{R}^3$. For local orientation preservation, we require

$$W \nearrow \infty \text{ as } J \searrow 0. \tag{2.2}$$

We assume material objectivity:

$$W(\mathbf{QF}, \mathbf{QD}, \mathbf{Qd}) \equiv W(\mathbf{F}, \mathbf{D}, \mathbf{d}) \ \ \forall \mathbf{Q} \in SO(3). \tag{2.3}$$

The total internal energy is given by $\int_\Omega W(\nabla\mathbf{f}, \nabla\mathbf{d}, \mathbf{d})dx$; the formal first-variation condition yields the following equilibrium equations

$$\begin{aligned} &div\left(W_F\right) + \mathbf{b}_o = \mathbf{0}, \\ &div\left(W_D\right) - W_d + \mathbf{m}_o = \mathbf{0}, \text{ in } \Omega, \end{aligned} \tag{2.4}$$

where $\mathbf{b}_o$, $\mathbf{m}_o$ denote prescribed body-force and body-couple fields acting on the system. The latter are readily incorporated in (2.4) via the principle of virtual work, for instance. For a given smooth tensor field $\mathbf{T} \in L(\mathbb{R}^2, \mathbb{R}^3)$, the divergence, denoted, $div\mathbf{T}$, is defined as the unique vector field satisfying $(div\mathbf{T}) \cdot \mathbf{a} = div(\mathbf{T}^T \mathbf{a})$ for all $\mathbf{a} \in \mathbb{R}^3$. Writing $\mathbf{T} = T_{i\gamma} \mathbf{e}_i \otimes \mathbf{e}_\gamma$, this implies $div(\mathbf{T}) = \frac{\partial T_{j\gamma}}{\partial x_\gamma} \mathbf{e}_j$. We define $\mathbf{N} := W_F = \frac{\partial W}{\partial F_{j\gamma}} \mathbf{e}_j \otimes \mathbf{e}_\gamma$, $\mathbf{M} := W_D = \frac{\partial W}{\partial D_{j\gamma}} \mathbf{e}_j \otimes \mathbf{e}_\gamma$, and $\boldsymbol{\mu} := W_d = \frac{\partial W}{\partial d_i} \mathbf{e}_i$; $\mathbf{N}$ and $\mathbf{M}$ are called the contact force and director-couple tensors (per unit undeformed length), respectively, while $\boldsymbol{\mu}$ is termed the intrinsic director couple (per unit undeformed area), cf. [20]. The two equations of (2.4) represent the local forms of force-balance and director-moment balance, respectively.

We assume a general class of "live" loadings of the form

$$\mathbf{b}_o = \hat{\mathbf{b}}(\lambda, \mathbf{F}, \mathbf{D}, \mathbf{d}, \mathbf{x}); \ \mathbf{m}_o = \hat{\mathbf{m}}(\lambda, \mathbf{F}, \mathbf{D}, \mathbf{d}, \mathbf{x}), \tag{2.5}$$

such that

$$\hat{\mathbf{b}}(0, \cdot) = \hat{\mathbf{m}}(0, \cdot) \equiv \mathbf{0}, \tag{2.6}$$

where $\lambda \in \mathbb{R}$ is a loading parameter. For simplicity, we assign Dirichlet b.c.'s:

$$\mathbf{f}\,|_{\partial\Omega} = \mathbf{x}, \ \mathbf{d}\,|_{\partial\Omega} = \mathbf{e}_3. \tag{2.7}$$

The flat reference configuration is presumed stress free:

$$W_F(\mathbf{I}_o, \mathbf{O}, \mathbf{e}_3) = W_D(\mathbf{I}_o, \mathbf{O}, \mathbf{e}_3) = \mathbf{O}, \ W_d(\mathbf{I}_o, \mathbf{O}, \mathbf{e}_3) = \mathbf{0}, \tag{2.8}$$

where $\mathbf{I}_o \in L(\mathbb{R}^2, \mathbb{R}^3)$ is the identity or inclusion map (cf. [9]) and $\mathbf{O} \in L(\mathbb{R}^2, \mathbb{R}^3)$ denotes zero.

For convenience, we define

$$\mathbf{P} = (\mathbf{F}, \mathbf{D}) \in L(\mathbb{R}^2, \mathbb{R}^3) \times L(\mathbb{R}^2, \mathbb{R}^3) \cong L(\mathbb{R}^2, \mathbb{R}^6), \tag{2.9}$$

and we sometimes write $W(\mathbf{P}, \mathbf{d}) \equiv W(\mathbf{F}, \mathbf{D}, \mathbf{d})$. With this in hand, strong ellipticity reads

$$\boldsymbol{\xi} \otimes \boldsymbol{\omega} \cdot \left(W_{PP}(\mathbf{P}, \mathbf{d})[\boldsymbol{\xi} \otimes \boldsymbol{\omega}]\right) > 0 \ \text{ for all } \boldsymbol{\xi} \in \mathbb{R}^6 \setminus \{\mathbf{0}\}, \ \boldsymbol{\omega} \in \mathbb{R}^2 \setminus \{\mathbf{0}\}, \tag{2.10}$$

which we assume holds throughout this work. In matrix form,

$W_{PP}(\mathbf{P}, \mathbf{d}) = \begin{bmatrix} W_{FF}(\mathbf{F}, \mathbf{D}, \mathbf{d}) & W_{FD}(\mathbf{F}, \mathbf{D}, \mathbf{d}) \\ W_{FD}(\mathbf{F}, \mathbf{D}, \mathbf{d}) & W_{DD}(\mathbf{F}, \mathbf{D}, \mathbf{d}) \end{bmatrix}$, and (2.10) is equivalent to

$$\begin{aligned} &\mathbf{a} \otimes \boldsymbol{\omega} \cdot \left(W_{FF}(\mathbf{F}, \mathbf{D}, \mathbf{d})[\mathbf{a} \otimes \boldsymbol{\omega}]\right) + 2\mathbf{a} \otimes \boldsymbol{\omega} \cdot \left(W_{FD}(\mathbf{F}, \mathbf{D}, \mathbf{d})[\mathbf{b} \otimes \boldsymbol{\omega}]\right) \\ &+ \mathbf{b} \otimes \boldsymbol{\omega} \cdot \left(W_{DD}(\mathbf{F}, \mathbf{D}, \mathbf{d})[\mathbf{b} \otimes \boldsymbol{\omega}]\right) > 0 \ \text{ for all } \mathbf{a}, \mathbf{b} \neq \mathbf{0} \in \mathbb{R}^3, \ \boldsymbol{\omega} \neq \mathbf{0} \in \mathbb{R}^2, \end{aligned} \tag{2.11}$$

where $W_{FF}(\cdot)[\mathbf{B}] = \frac{\partial^2 W(\cdot)}{\partial F_{i\mu} \partial F_{j\nu}} B_{j\nu} \mathbf{e}_i \otimes \mathbf{e}_\mu$, $W_{FD}(\cdot)[\mathbf{B}] = \frac{\partial^2 W(\cdot)}{\partial F_{i\mu} \partial D_{j\nu}} B_{j\nu} \mathbf{e}_i \otimes \mathbf{e}_\mu$, etc., and $\mathbf{A} \cdot \mathbf{B} := tr(\mathbf{A}^T \mathbf{B})$ for all $\mathbf{A}, \mathbf{B} \in L(\mathbb{R}^2, \mathbb{R}^3)$.

At the reference configuration we require a stronger but physically reasonable condition: For convenience, define $W_{FF}^o := W_{FF}(\mathbf{I}_o, \mathbf{O}, \mathbf{e}_3)$, etc., and consider the matrix

$$\begin{bmatrix} W_{FF}^o & W_{FD}^o & W_{Fd}^o \\ W_{FD}^o & W_{DD}^o & W_{Dd}^o \\ W_{FD}^o & W_{Dd}^o & W_{dd}^o \end{bmatrix}, \tag{2.12}$$

which defines a linear transformation, denoted $\mathrm{K}_o$, of the vector space

$$Y = L(\mathbb{R}^2, \mathbb{R}^3) \times L(\mathbb{R}^2, \mathbb{R}^3) \times \mathbb{R}^3, \tag{2.13}$$

into itself, i.e., $\mathrm{K}_o \in L(Y)$. As shown in the Appendix, (2.3) implies that the null space of $\mathrm{K}_o$ contains the subspace

$$\mathrm{N} = \{(\mathbf{W}\mathbf{I}_o, \mathbf{O}, \mathbf{W}\mathbf{e}_3) \in Y : \mathbf{W} \in Skew(\mathbb{R}^3)\}, \tag{2.14}$$

where $Skew(\mathbb{R}^3) = \{\mathbf{W} \in L(\mathbb{R}^3) : \mathbf{W}^T = -\mathbf{W}\}$. Denoting a typical element $\xi \in Y$ via $\boldsymbol{\xi} = (\mathbf{A}, \mathbf{B}, \mathbf{a})$, we assume that the restriction of $\mathrm{K}_o$ to $\mathrm{N}^\perp$ is positive definite:

$$\langle \boldsymbol{\xi}, \mathrm{K}_o \boldsymbol{\xi} \rangle > 0 \ \forall \ \boldsymbol{\xi} \in \mathrm{N}^\perp \setminus \{\mathbf{0}\}, \tag{2.15}$$

where for any $\boldsymbol{\xi} = (\mathbf{A}, \mathbf{B}, \mathbf{a})$, $\boldsymbol{\zeta} = (\mathbf{C}, \mathbf{D}, \mathbf{b}) \in Y$, we define the inner product $\langle \boldsymbol{\xi}, \boldsymbol{\zeta} \rangle := \mathbf{A} \cdot \mathbf{C} + \mathbf{B} \cdot \mathbf{D} + \mathbf{a} \cdot \mathbf{b}$. Accordingly,

$$\begin{aligned} \mathrm{N}^\perp &= \{(\mathbf{A}, \mathbf{B}, \mathbf{a}) \in Y : \mathbf{A} \cdot (\mathbf{W}\mathbf{I}_o) + \mathbf{a} \cdot (\mathbf{W}\mathbf{e}_3) = 0 \ \forall \ \mathbf{W} \in Skew(\mathbb{R}^3)\} \\ &= \{(\mathbf{H}\mathbf{I}_o, \mathbf{B}, \mathbf{H}\mathbf{e}_3) : \mathbf{H} \in Sym(\mathbb{R}^3), \mathbf{B} \in L(\mathbb{R}^2, \mathbb{R}^3)\}, \end{aligned} \tag{2.16}$$

where $Sym(\mathbb{R}^3) = \{\mathbf{H} \in L(\mathbb{R}^3) : \mathbf{H}^T = \mathbf{H}\}$.

**Remark 2.1.** With (2.15) in hand, we see that $\mathrm{N}$ is precisely $Null(\mathrm{K}_o)$. Assumption (2.15) is analogous to the requirement in 3D elasticity that the elasticity tensor at a stress-free reference configuration is positive definite on symmetric tensors, e.g., [8].

Before expressing the field equations in terms of a nonlinear operator, we introduce

$$\mathbf{f}(\mathbf{x}) = \mathbf{x} + \mathbf{u}^1(\mathbf{x}); \ \ \mathbf{d}(\mathbf{x}) = \mathbf{e}_3 + \mathbf{u}^2(\mathbf{x}). \tag{2.17}$$

Then (2.4), (2.5) and (2.7) lead to the following boundary value problem:

$$\begin{aligned} &div\{W_F(\mathbf{I}_o + \nabla \mathbf{h}^1, \nabla \mathbf{h}^2, \mathbf{e}_3 + \mathbf{h}^2)\} + \mathbf{b}(\lambda, \nabla \mathbf{u}^1, \nabla \mathbf{u}^2, \mathbf{u}^2, \mathbf{x}) = \mathbf{0}, \\ &div\{W_D(\mathbf{I}_o + \nabla \mathbf{h}^1, \nabla \mathbf{h}^2, \mathbf{e}_3 + \mathbf{h}^2)\} - W_d(\mathbf{I}_o + \nabla \mathbf{h}^1, \nabla \mathbf{h}^2, \mathbf{e}_3 + \mathbf{h}^2) \\ &\qquad + \mathbf{m}(\lambda, \nabla \mathbf{u}^1, \nabla \mathbf{u}^2, \mathbf{u}^2, \mathbf{x}) = \mathbf{0}, \ \text{ in } \Omega, \\ &\qquad \mathbf{u}^1|_{\partial\Omega} = \mathbf{u}^2|_{\partial\Omega} = \mathbf{0}, \end{aligned} \tag{2.18}$$

where

$$\begin{aligned}&\mathbf{b}(\lambda,\nabla\mathbf{u}^1,\nabla\mathbf{u}^2,\mathbf{u}^2,\mathbf{x}):=\widehat{\mathbf{b}}(\lambda,\mathbf{I}_o+\nabla\mathbf{u}^1,\nabla\mathbf{u}^2,\mathbf{e}_3+\mathbf{u}^2,\mathbf{x}),\ \text{and}\\&\mathbf{m}(\lambda,\nabla\mathbf{u}^1,\nabla\mathbf{u}^2,\mathbf{u}^2,\mathbf{x}):=\widehat{\mathbf{m}}(\lambda,\mathbf{I}_o+\nabla\mathbf{u}^1,\nabla\mathbf{u}^2,\mathbf{e}_3+\mathbf{u}^2,\mathbf{x}).\end{aligned}\tag{2.19}$$

We then write (2.18) in quasi-linear form via

$$\begin{aligned}\mathcal{C}(\nabla\mathbf{u},\mathbf{u}^2)\{\nabla^2\mathbf{u}\}+\Upsilon(\lambda,\nabla\mathbf{u},\mathbf{u}^2)=\mathbf{0}\quad&\text{in }\Omega,\\ \mathbf{u}|_{\partial\Omega}=\mathbf{0},&\end{aligned}\tag{2.20}$$

where $\mathbf{u}:=(\mathbf{u}^1,\mathbf{u}^2)$, $\nabla\mathbf{u}:=(\nabla\mathbf{u}^1,\nabla\mathbf{u}^2)$, $\nabla^2\mathbf{u}:=(\nabla^2\mathbf{u}^1,\nabla^2\mathbf{u}^2)$, and $\mathcal{C}(\nabla\mathbf{u},\mathbf{u}^2)\{\nabla^2\mathbf{u}\}:=$ $W_{PP}(\mathbf{I}_o+\nabla\mathbf{u}^1,\nabla\mathbf{u}^2,\mathbf{e}_3+\mathbf{u}^2)\{\nabla^2\mathbf{u}\}$. In matrix form, this reads

$$\mathcal{C}(\nabla\mathbf{u},\mathbf{u}^2)\{\nabla^2\mathbf{u}\}=\begin{bmatrix}W_{FF}(\mathbf{I}_o+\nabla\mathbf{u}^1,\nabla\mathbf{u}^2,\mathbf{e}_3+\mathbf{u}^2)+W_{FD}(\mathbf{I}_o+\nabla\mathbf{u}^1,\nabla\mathbf{u}^2,\mathbf{e}_3+\mathbf{u}^2)\\W_{FD}(\mathbf{I}_o+\nabla\mathbf{u}^1,\nabla\mathbf{u}^2,\mathbf{e}_3+\mathbf{u}^2)+W_{DD}(\mathbf{I}_o+\nabla\mathbf{u}^1,\nabla\mathbf{u}^2,\mathbf{e}_3+\mathbf{u}^2)\end{bmatrix}\begin{pmatrix}\{\nabla^2\mathbf{u}^1\}\\\{\nabla^2\mathbf{u}^2\}\end{pmatrix},\tag{2.21}$$

where $W_{FF}(\cdot)\{\nabla^2\mathbf{u}^1\}:=\dfrac{\partial^2W(\cdot)}{\partial F_{i\mu}\partial F_{j\nu}}\dfrac{\partial^2u_j^1}{\partial x_\mu\partial x_\nu}\mathbf{e}_i$, etc. The second term on the left side of $(2.20)_1$ captures the lower-order terms:

$$\begin{aligned}\Upsilon(\lambda,\nabla\mathbf{u},\mathbf{u}^2):=[&W_{Fd}(\mathbf{I}_o+\nabla\mathbf{u}^1,\nabla\mathbf{u}^2,\mathbf{e}_3+\mathbf{u}^2)[\nabla\mathbf{u}^2]+\mathbf{b}(\lambda,\nabla\mathbf{u}^1,\nabla\mathbf{u}^2,\mathbf{u}^2),\\&W_{Dd}(\mathbf{I}_o+\nabla\mathbf{u}^1,\nabla\mathbf{u}^2,\mathbf{e}_3+\mathbf{u}^2)[\nabla\mathbf{u}^2]-W_d(\mathbf{I}_o+\nabla\mathbf{u}^1,\nabla\mathbf{u}^2,\mathbf{e}_3+\mathbf{u}^2)+\mathbf{m}(\lambda,\nabla\mathbf{u}^1,\nabla\mathbf{u}^2,\mathbf{u}^2)]\end{aligned}\tag{2.22}$$

where $W_{Dd}(\cdot)[\nabla\mathbf{u}^2]:=\dfrac{\partial^2W(\cdot)}{\partial F_{i\mu}d_j}\dfrac{\partial u_j^2}{\partial x_\mu}\mathbf{e}_i$, etc.

We assume

$$\begin{aligned}&W\in C^4(\Theta,\mathbb{R});\\&\mathbf{b},\mathbf{m}\in C^3(\mathbb{R}\times L(\mathbb{R}^2,\mathbb{R}^3)\times L(\mathbb{R}^2,\mathbb{R}^3)\times\mathbb{R}^3\times\mathbb{R}^3,\mathbb{R}^3).\end{aligned}\tag{2.23}$$

Finally, let $C^{k,\alpha}(\bar{\Omega},\mathbb{R}^6)$ denote the usual Hölder space with norm $\|\ \|_{k,\alpha}$, for $k=0,1,2,...,\ 0<\alpha<1$, cf. [7]. We then define the Banach spaces

$$\begin{aligned}&V:=C^\alpha(\bar{\Omega},\mathbb{R}^6),\ \|\ \|_\alpha\ (k=0);\\&U:=\{\mathbf{u}\in C^{2,\alpha}(\bar{\Omega},\mathbb{R}^6):\mathbf{u}|_{\partial\Omega}=\mathbf{0}\},\ \|\ \|_{2,\alpha},\end{aligned}\tag{2.24}$$

and the open subset

$$\mathcal{E}=\{\mathbf{u}\in U:J=\det[\mathbf{I}+(\nabla\mathbf{u}^1+\mathbf{u}^2\otimes\mathbf{e}_3)]>0\ \text{ on }\bar{\Omega}\},\tag{2.25}$$

where $\mathbf{I}=\mathbf{I}_o+\mathbf{e}_3\otimes\mathbf{e}_3$ is the identity on $\mathbb{R}^3$. The left side of (2.20) then defines a mapping $F:\mathbb{R}\times\mathcal{E}\to V$, so that (2.20) (or (2.18)) is expressed abstractly as

$$F(\lambda,\mathbf{u})=0.\tag{2.26}$$

**3. Local Solution Path**

In view of (2.6), (2.18) and (2.19), we see that $(\lambda, \mathbf{u}) = (0, \mathbf{0}) \in \mathcal{E}$ is a solution point of (2.26), i.e.,

$$F(0, \mathbf{0}) = 0. \tag{3.1}$$

Our goal in this section is to prove the existence of a local solution path containing $(0, \mathbf{0})$ via the implicit function theorem. The presumed smoothness (2.23) ensures that $F(\cdot)$ is Fréchet differentiable, cf. [5], [22]. Thus, we may obtain the Fréchet derivative of $\mathbf{u} \mapsto F$, denoted $D_u F$, by taking an arbitrary directional or Gateaux derivative at $(\lambda, \mathbf{u}) = (0, \mathbf{0})$, viz.,

$$T_o[\mathbf{h}] := D_u F(0, \mathbf{0})[\mathbf{h}] = \frac{d}{d\varepsilon} F(0, \varepsilon \mathbf{h})\,|_{\varepsilon=0},$$

for all $\mathbf{h} = (\mathbf{h}^1, \mathbf{h}^2) \in U$, where $T_o \in L(U, V)$. We employ (2.6), (2.8), (2.18), (2.19) and (2.22) to find

$$T_o[\mathbf{h}] = \begin{pmatrix} div\{W^o_{FF}[\nabla \mathbf{h}^1] + W^o_{FD}[\nabla \mathbf{h}^2] + W^o_{Fd}[\mathbf{h}^2]\}, \\ div\{W^o_{FD}[\nabla \mathbf{h}^1] + W^o_{DD}[\nabla \mathbf{h}^2] + W^o_{Dd}[\mathbf{h}^2]\} - \{W^o_{FD}[\nabla \mathbf{h}^1] - W^o_{Dd}[\nabla \mathbf{h}^2] + W^o_{dd}[\mathbf{h}^2]\} \end{pmatrix}, \tag{3.2}$$

**Lemma 3.1.** *Assuming the hypotheses of Section* 2**,** $T_o \in L(U, V)$ *is injective*.

*Proof.* Consider $T_o[\mathbf{h}] = 0.$ We take the dot product of $T_o[\mathbf{h}]$ with $\mathbf{h}$, viz., we dot the first component of (3.2) with $\mathbf{h}^1$ and add the result to the dot product of the second component with $\mathbf{h}^2$; we integrate the resulting sum over $\Omega$. Integration by parts reveals

$$\int_\Omega \langle \boldsymbol{\psi}, \mathrm{K}_o \boldsymbol{\psi} \rangle dx = 0, \tag{3.3}$$

where $\boldsymbol{\psi} := (\nabla \mathbf{h}^1, \nabla \mathbf{h}^2, \mathbf{h}^2)$, and $\mathrm{K}_o$ is the matrix (2.12). In view of hypothesis (2.15), equation (3.3) implies that the field $\boldsymbol{\psi}$ must take its values in $Null(\mathrm{K}_o) = \mathrm{N}$, cf. (2.14) and Remark 2.1. Hence, $\nabla \mathbf{h}^1 = \mathbf{W}\mathbf{I}_o, \nabla \mathbf{h}^2 = \mathbf{O}, \mathbf{h}^2 = \mathbf{W}\mathbf{e}_3$, for $\mathbf{W} : \bar{\Omega} \to Skew(\mathbb{R}^3)$. The second and third of these imply

$$\mathbf{h}^2 = \mathbf{W}\mathbf{e}_3 = \mathbf{c} \in \mathbb{R}^3 \Rightarrow \mathbf{e}_\gamma \cdot \mathbf{W}\mathbf{e}_3 = c_\gamma,\ \gamma = 1, 2. \tag{3.4}$$

The infinitesimal in-plane strain is given by $\mathbf{E} = \frac{1}{2}(\mathbf{I}_o^T \nabla \mathbf{h}^1 + [\nabla \mathbf{h}^1]^T \mathbf{I}_o)$, cf. [9]; this vanishes at $\mathbf{E} =$ $\nabla \mathbf{h}^1(\mathbf{x}) = \mathbf{W}(x)\mathbf{I}_o$ :

$$\mathbf{E} = \frac{1}{2}(\mathbf{I}_o^T \mathbf{W}\mathbf{I}_o + [\mathbf{W}\mathbf{I}_o]^T \mathbf{I}_o) = \frac{1}{2}(\mathbf{I}_o^T \mathbf{W}\mathbf{I}_o + \mathbf{I}_o^T \mathbf{W}^T \mathbf{I}_o) \equiv \boldsymbol{\Theta},$$

where $\boldsymbol{\Theta} \in L(\mathbb{R}^2)$ denotes zero. From the compatibility conditions of classical elasticity [26], we then conclude that $\mathbf{I}_o^T \mathbf{W}(\mathbf{x})\mathbf{I}_o = \boldsymbol{\Omega} \in Skew(\mathbb{R}^2)$. Combing this with (3.4), we see that that skew field $\mathbf{W}$ is constant, say, $\mathbf{W} = \mathbf{W}_o \in Skew(\mathbb{R}^3)$. Thus, $\mathbf{h}^1 = \mathbf{W}_o \mathbf{I}_o \mathbf{x} + \mathbf{c}_o$, $\mathbf{c}_o \in \mathbb{R}^3$, and $\mathbf{h}^2 = \mathbf{W}_o \mathbf{e}_3$. From the boundary conditions in (2.18), we now conclude that $\mathbf{h}^1 = \mathbf{h}^2 \equiv \mathbf{0}$. □

**Lemma 3.2.** *Given the hypotheses of Section* 2, $T_o \in L(U,V)$ *is surjective*.

*Proof.* Let $\mathcal{I}: L(\mathbb{R}^2,\mathbb{R}^3) \to L(\mathbb{R}^2,\mathbb{R}^3)$ denote the identity. Then $(\mathcal{I},\mathcal{I})$ is the identity on $L(\mathbb{R}^2,\mathbb{R}^3) \times L(\mathbb{R}^2,\mathbb{R}^3)$. Now consider the one-parameter family of fourth-order tensors

$$\mathcal{C}_\mu := \mu(\mathcal{I},\mathcal{I}) + (1-\mu)\mathcal{C}_o, \;\; \mu \in [0,1], \tag{3.5}$$

where $\mathcal{C}_o := \mathcal{C}(\mathbf{O},\mathbf{0})$, cf. (2.21). Clearly, $\mathcal{C}_\mu$ satisfies strong ellipticity as in (2.10), uniformly in $\mu \in [0,1]$. Next, define a linear operator $T_\mu \in L(U,V)$ via the convex combination

$$T_\mu[\mathbf{h}] := div\mathcal{C}_\mu[\nabla\mathbf{h}] = \mu\Delta\mathbf{h} + (1-\mu)\{div\mathcal{C}_o[\nabla\mathbf{h}] + ...\}, \tag{3.6}$$

where the complete term multiplying $(1-\mu)$ is given in (3.2); only the principal part is shown above. Also, $\Delta\mathbf{h} := (\Delta\mathbf{h}^1, \Delta\mathbf{h}^2) = (div\nabla\mathbf{h}^1, div\nabla\mathbf{h}^2)$. At $\mu = 1$, we have $T_1[\mathbf{h}] = \Delta\mathbf{h}$, and we consider the boundary value problem

$$\begin{aligned} &\Delta\mathbf{h} = \mathbf{g} \text{ in } \Omega, \\ &\mathbf{h}\,|_{\partial\Omega} = \mathbf{0}. \end{aligned} \tag{3.7}$$

This represents six uncoupled Poisson problems. With $\mathbf{g} \in V$, each such component problem has a unique solution in $C^{2,\alpha}(\bar{\Omega}) \cap \{h\,|_{\partial\Omega} = 0\}$, cf. [7]. Hence, (3.7) has a unique solution $\mathbf{h} \in U$. In particular, $T_1$ has Fredholm index zero; the Fredholm index of $T_\mu$ is defined as $\dim Null(T_\mu) -$ codim$Range(T_\mu)$.

Next, the Schauder estimates for (3.6) read:

$$\|\mathbf{h}\|_{2,\alpha} \le C\left[\|T_u[\mathbf{h}]\|_\alpha + \|\mathbf{h}\|_\infty\right] \text{ for all } \mathbf{h} \in U, \tag{3.8}$$

where $\|\;\|_\infty$ denotes the maximum norm, and the constant $C$ is independent of $\mathbf{h}$ and $\mu$, cf. [1]. Inequality (3.8) ensures that the Fredholm index is constant on $\mu \in [0,1]$, cf. [16]. By Lemma 3.1, the null space of $T_o$ is trivial. Thus, codim$Range(T_o) = 0$ □

With Lemmas 3.1 and 3.2 in hand, the implicit function theorem yields:

**Theorem 3.1.** *There exists a unique local path of solutions of* (2.26) (*or* (2.18)) *of the form*

$$\mathcal{P}_\varepsilon := \{(\lambda, \tilde{\mathbf{u}}(\lambda)) : |\lambda| < \varepsilon\}, \;\; F(\lambda, \tilde{\mathbf{u}}(\lambda)) \equiv 0, \tag{3.9}$$

*where* $\varepsilon > 0$ *is sufficiently small,* $\lambda \mapsto \tilde{\mathbf{u}}(\lambda)$ *is* $C^1$, *and* $\tilde{\mathbf{u}}(0) = \mathbf{0}$. *Moreover,* (3.9) *captures all solutions of* (2.26) *in a sufficiently small neighborhood of* $(0,\mathbf{0}) \in \mathbb{R} \times \mathcal{E}$.

*Proof.* Given that $T_o \in L(U,V)$ is invertible, it follows directly from the implicit function theorem that that $\tilde{\mathbf{u}}(\lambda) \in \mathcal{U}$ on $(-\varepsilon,\varepsilon)$. It only remains to show that $\tilde{\mathbf{u}}(\lambda) \in \mathcal{E}$ for all $|\lambda| < \varepsilon$. Consider the expression for $J$ in (2.25), evaluated along $\tilde{\mathbf{u}}(\lambda)$, which we denote by $\tilde{J}_\lambda(\mathbf{x}) = \det[\mathbf{I} + \tilde{\mathbf{\Pi}}_\lambda(\mathbf{x})]$, where

$$\mathbf{\Pi} := (\nabla \mathbf{u}^1 + \mathbf{u}^2 \otimes \mathbf{e}_3). \tag{3.10}$$

Since $C^{1,\alpha}(\bar{\Omega},\mathbb{R}^3)$ is a Banach algebra, we see that $\tilde{J}_\lambda \in C^{1,\alpha}(\bar{\Omega})$ on $(-\varepsilon,\varepsilon)$. For sufficiently small $\varepsilon > 0$, $\|\tilde{\mathbf{u}}(\lambda)\|_{2,\alpha}$ can be made small enough so that $\max_{x\in\bar{\Omega}} |\tilde{\mathbf{\Pi}}_\lambda(\mathbf{x})| < 1$, which, in turn, ensures that $\tilde{J}_\lambda > 0$ on $\bar{\Omega}$, $|\lambda| < \varepsilon$. □

**4. Global Solution Branch**

Next, we show that the local solution path (3.9) is part of a global branch of solutions, viz., a connected, locally compact set of solution pairs of (2.26). To begin, we define $\mathcal{O}$ as the maximal connected set in $\mathcal{E}$ containing $\mathbf{u} = \mathbf{0}$, called the component of $\mathbf{0}$ in $\mathcal{E}$, denoted

$$\mathcal{O} := comp\{\mathbf{0}\} \text{ in } \mathcal{E}. \tag{4.1}$$

In addition, for each $\delta > 0$, we define (using the notation (3.10))

$$\mathcal{O}_\delta := \{\mathbf{u} \in \mathcal{O} : J = \det[\mathbf{I} + \mathbf{\Pi}] > \delta\}. \tag{4.2}$$

Note that each $\mathcal{O}_\delta \subset U$ is open, $\overline{\mathcal{O}_\delta} \subset \mathcal{O}$, and $\mathcal{O} = \bigcup_{\delta>0} \mathcal{O}_\delta$. Hence, $\mathcal{O} \subset U$ is open.

**Remark 4.1** The proof of Theorem 3.1 implies that the local solution path satisfies

$$\mathcal{P}_\varepsilon \subset \mathbb{R} \times \mathcal{O},$$

for sufficiently small $\varepsilon > 0$.

Let $\mathcal{S}$ denote the solution set of (26) in $\mathbb{R} \times \mathcal{E}$ and define

$$\Sigma := comp\{(0,\mathbf{0})\} \text{ in } \mathcal{S}. \tag{4.3}$$

We now state a global theorem:

**Theorem 4.1.** *Given the hypotheses of Section* 2, *the solution branch* $\Sigma$ *is characterized by at least one of the following:*

(i) $\Sigma$ *is unbounded in* $\mathbb{R} \times (U \cap C^2(\overline{\Omega},\mathbb{R}^6)$.

(ii) $\Sigma \setminus \{(0,\mathbf{0})\}$ *is connected.*

(iii) $\Sigma \not\subset \mathbb{R} \times \mathcal{O}_\delta$ *for each* $\delta > 0$.

We note that if alternative (ii) holds, then there at least one solution point $(0,\mathbf{u}) \in \mathcal{C}$ with $\mathbf{u} \neq \mathbf{0}$, i.e., the unloaded state admits a nontrivial solution, cf. (2.4)-(2.7). If (iii) is true, then there is a sequence of solution points

$$\{(\lambda_j, \mathbf{u}_j)\} \subset \Sigma \text{ such that } \inf_{x\in\bar{\Omega}} J_j \searrow 0, \tag{4.4}$$

where $J_j := \det[\mathbf{I} + \mathbf{\Pi}_j]$, cf. (3.10). If $\{\mathbf{u}_j\}$ is bounded in this case, then there is a subsequence that converges to, say, $\mathbf{u}_*$ in the $C^2$ topology. Moreover, $\{J_j\}$ is bounded in $C^{1,\alpha}(\bar{\Omega})$, implying that $J_j \to J_*$ in the usual $C^1$ topology. Hence, from (4.4) we obtain

$$J_*(\mathbf{x}_*) = \det[\mathbf{I} + \mathbf{\Pi}(\mathbf{x}_*)] = 0, \tag{4.5}$$

for one or more points $\mathbf{x}_* \in \bar{\Omega}$. In this case, we also note from (2.7) that $J_j \equiv 3$ on $\partial\Omega$ for each $j$, implying

$$\|J_*\|_\infty \geq 1. \tag{4.6}$$

To prove the theorem, we need to first demonstrate two crucial properties of the mapping $F$ in (2.26), which we provide in the following Propositions. To begin, we denote the Fréchet derivative of $\mathbf{u} \mapsto F$ at $(\lambda, \mathbf{u}) \in \mathbb{R} \times \mathcal{O}$ by $T(\lambda, \mathbf{u}) := D_u F(\lambda, \mathbf{u}) \in L(U, V)$. From (2.20) and (2.23), this takes the form

$$T(\lambda, \mathbf{u})[\mathbf{h}] = \mathcal{C}(\nabla\mathbf{u}, \mathbf{u}^2)\{\nabla^2\mathbf{h}\} + ..., \tag{4.7}$$

for all $\mathbf{h} \in U$, where only the principal part of the linear operator is shown.

**Proposition 4.1.** *For each* $(\lambda, \mathbf{u}) \in \mathbb{R} \times \mathcal{O}$, $T(\lambda, \mathbf{u}) \in L(U, V)$ *is a Fredholm operator of index zero*.

*Proof.* The argument is like that used in the proof of Lemma 3.2. For fixed $(\lambda, \mathbf{u}) \in \mathbb{R} \times \mathcal{O}$, the Schauder estimates for (4.7) read

$$\|\mathbf{h}\|_{2,\alpha} \leq C\left[\|T(\lambda, \mathbf{u})[\mathbf{h}]\|_\alpha + \|\mathbf{h}\|_\infty\right], \tag{4.8}$$

for all $\mathbf{h} \in U$, where the constant $C$ is independent of $\mathbf{h}$, cf. [1]. Inequality (4.8) holds if and only if $T(\lambda, \mathbf{u}) \in L(U, V)$ is semi-Fredholm, viz., the null space is finite-dimensional and the range is closed, cf. [27]. Since $\mathbb{R} \times \mathcal{O}$ is path connected, the assertion now follows from the observation that $T(0, \mathbf{0}) = T_o$ : The continuity of the Fredholm index implies that the index of $T(\lambda, \mathbf{u})$ equals that of $T_o$, which is zero courtesy of Lemmas 3.1 and 3.2 □

**Proposition 4.2** *For each* $\delta > 0$, $F : \mathbb{R} \times \mathcal{O}_\delta \to V$ *is proper, i.e.,* $F^{-1}(\Upsilon) \cap \bar{B}$ *is compact for each compact set* $\Upsilon \subset V$ *and bounded set* $B \subset \mathbb{R} \times \overline{\mathcal{O}_\delta}$.

*Proof.* We use an argument like that employed in [11]. We first rewrite $F$ in operator form, isolating the quasilinear part from the nonlinear lower-order terms:

$$F(\lambda,u)=\mathcal{L}(\mathbf{u})[\mathbf{u}]+\mathcal{G}(\lambda,\mathbf{u}),$$

where for fixed $\mathbf{u}\in U$, $\mathcal{L}(\mathbf{u})[\mathbf{h}]:=\mathcal{C}(\nabla\mathbf{u},\mathbf{u}^2)\{\nabla^2\mathbf{h}\}$ for all $\mathbf{h}\in U$, and $\mathcal{G}(\lambda,\mathbf{u}):=\Upsilon(\lambda,\nabla\mathbf{u},\mathbf{u}^2)$, cf. (2.22). Set $F(\lambda_j,\mathbf{u}_j)=p_j$, where $\{p_j\}\subset V$ converges, say, $p_j\to p_*$, and $\{(\lambda_j,\mathbf{u}_j)\}\subset B\subset\mathbb{R}\times\overline{\mathcal{O}_\delta}$ is bounded. Our goal is to show that $\{(\lambda_j,\mathbf{u}_j)\}$ has a convergent subsequence. By compact embedding, we deduce for possibly a subsequence (not relabeled) that

$$\lambda_j\to\lambda_*,\ \text{and}\ \mathbf{u}_j\to\mathbf{u}_*\ \text{in}\ C^{1,\alpha}(\bar{\Omega},\mathbb{R}^6),$$

implying

$$\mathcal{G}(\lambda_j,\mathbf{u}_j)\to\mathcal{G}(\lambda_*,\mathbf{u}_*)\ \text{in}\ V;$$
$$\left\|[\mathcal{L}(\mathbf{u}_j)-\mathcal{L}(\mathbf{u}_*)][\mathbf{u}_j]\right\|_\alpha\le\varepsilon\left\|\mathbf{u}_j\right\|_{2,\alpha}\le\varepsilon M\quad\text{for all}\ j\ge N_\varepsilon,$$

for some sufficiently large number $N_\varepsilon$. Since $\{\mathbf{u}_j\}$ is bounded in $U$, the coefficients of the quasilinear operator are uniformly continuous in the index $j$. Thus, the above estimate implies

$$\mathcal{L}(\mathbf{u}_j)[\mathbf{u}_j]-\mathcal{L}(\mathbf{u}_*)[\mathbf{u}_j]\to 0\ \text{in}\ V.$$

In view of (2.10) and (2.21), it follows that $\mathcal{L}(\mathbf{u}_*)\in L(U,V)$ defines a uniformly elliptic, formally self-adjoint, sectorial operator. Let $\xi\in\mathbb{R}$ belong to the resolvent set, i.e., the shifted operator $\mathcal{L}(\mathbf{u}_*)-\xi\in L(U,V)$ is invertible. Finally, consider

$$\{\mathcal{L}(\mathbf{u}_*)-\xi\}[\mathbf{u}_j]=\mathcal{L}(\mathbf{u}_*)[\mathbf{u}_j]-\mathcal{L}(\mathbf{u}_j)[\mathbf{u}_j]-\mathcal{G}(\lambda_j,\mathbf{u}_j)+p_j-\xi\mathbf{u}_j,$$

which converges to $-\mathcal{G}(\lambda_*,\mathbf{u}_*)+p_*-\xi\mathbf{u}_*$ in $V$, i.e., $\{\mathcal{L}(\mathbf{u}_*)-\xi\}[\mathbf{u}_j]$ converges in $V$. Since $\mathcal{L}(\mathbf{u}_*)-\xi$ is a bijection, this is equivalent to the convergence of $\{\mathbf{u}_j\}$ in $\overline{\mathcal{O}_\delta}\subset U$. □

With Propositions 4.1 and 4.2 in hand, we may employ the $C^1$ nonlinear Fredholm degree developed first in the context of $C^2$ maps in [6], later generalized to continuously differentiable, nonlinear Fredholm maps in [24]. This base-point degree, henceforth called the FPR degree, possesses many of the usual properties of the Leray-Schauder degree. We give a brief description: Given, say, $G:\mathcal{D}\subset U\to V$, with $\mathcal{D}$ open and bounded, we consider the equation $G(u)=0,\ u\in\mathcal{D}$. Assume that $G$ possesses the properties established in Propositions 4.1 and 4.2, viz., $G|_{\bar{\mathcal{D}}}$ is proper and $DG(u)\in L(U,V)$ is Fredholm of index zero. We denote the FPR degree of $G$ with respect to $\mathcal{D}$ and 0 via $\deg_p(G,\mathcal{D}):=\deg_p(G,\mathcal{D},0)$, where $p\in\mathcal{D}$ is called a base point. In the regular-value case, a base point $p\in\mathcal{D}$ is chosen to be path connected to each of the finitely many solutions of $G(u)=0,\ u\in\mathcal{D}$. If no base point exists and/or if $\mathcal{D} \cap G^{-1}(0) = \varnothing$, then $\deg_p(G, \mathcal{D}) := 0$. It can be shown that a change in base point delivers the same integer value to within a sign, cf. [6]. As such, the FPR degree is not homotopy invariant. For our purposes here, it's enough to note that its absolute value $|\deg_p(G, \mathcal{D})|$ is homotopy invariant, and moreover, that $\deg_p(G, \mathcal{D}) \neq 0 \Rightarrow G(u) = 0$ has at least one solution $u \in \mathcal{D}$.

A detailed general proof of Theorem 4.1 without the complication of alternative (iii) is provided in both [17] and in [10]. The former presumes a general Banach-space setting, while employing an oriented $C^2$ Fredholm degree. The FPR degree is employed in [10] for problems of incompressible nonlinear elasticity, where property (iii) of Theorem 4.1 is not an issue $(J \equiv 1)$. Here, we use an argument like that given in [15] to account for (iii). A sketch of the proof is as follows:

Clearly the local solution from Theorem 3.1 is contained in $\Sigma$, i.e., $\mathcal{P}_\varepsilon \subset \Sigma$. We argue by contradiction: Suppose that $\Sigma$ is not characterized by any of the alternatives (i)-(iii). In particular, $\Sigma$ is bounded in $\mathbb{R} \times U$ with $\Sigma \setminus \{(0, \mathbf{0})\}$ not connected. Then $\Sigma$ is compact by properness, and by a well-known argument [25], we can find a bounded open set $\Xi \subset \mathbb{R} \times \mathcal{O}$ such that $\bar{\Xi} \cap \mathcal{S} = \Sigma$ and $\partial\Xi \cap \mathcal{S} = \varnothing$. Define $\Xi_\lambda := \{\mathbf{u} \in \mathcal{O} : (\lambda, \mathbf{u}) \in \Xi\}$. By homotopy invariance, it follows that $|\deg_0(F(\lambda, \cdot), \Xi_\lambda)| = const.$, where we have employed base point $p = \mathbf{0}$, as in [TJH]. Since $\Xi_\lambda = \varnothing$ for sufficiently large $|\lambda|$, we conclude that

$$|\deg_0(F(\lambda, \cdot), \Xi_\lambda)| \equiv 0. \qquad (4.9)$$

On the other hand, the uniqueness part of implicit function theorem implies that $\Xi_0 \cap F^{-1}(0) = \{(0, \mathbf{0})\}$, cf. Theorem 3.1. Consequently, $|\deg_0(F(0, \cdot), \Xi_0)| = 1$, which contradicts (4.9).

We are then left with at least one of the alternatives (i)-(iii) of Theorem 4.1. For case (i), the construction implies unboundedness in the $\mathbb{R} \times U$ topology. However, boundedness in $\mathbb{R} \times C^2(\overline{\Omega}, \mathbb{R}^6)$ contradicts presumed unboundedness in $\mathbb{R} \times U$, as shown in [14] (in the context of 3D nonlinear elasticity). The argument is the same here, which we do not repeat. Assuming that (iii) does not hold, the proof accounting for (ii) is more delicate than that given above. We refer to [17] and [10] for details. To realize alternative (iii), we repeat the same arguments above on $\mathbb{R} \times \mathcal{O}_\delta$ for each $\delta > 0$. As it stands, property (iii) implying (4.4) is possible with or without (i) and/or (ii) being true.

We finish with a special class of stored energies for which (4.4) is only possible when property (i) of Theorem 4.1 holds, viz., $\{(\lambda_j, \mathbf{u}_j)\} \subset \Sigma$ as given in (4.4) implies that $\{\mathbf{u}_j\}$ is unbounded. As a first step, we observe that the energy density is independent of $\mathbf{x} \in \overline{\Omega}$. By Noether's theorem (cf. [4], [23]), any classical solution of (2.4), (2.5) is also a classical solution of the balance law

$$div[W(\mathbf{P},\mathbf{d})\mathbf{1}-\mathbf{P}^TW_P(\mathbf{P},\mathbf{d})]=\mathbf{P}^T(\tilde{\mathbf{b}}(\lambda,\mathbf{P},\mathbf{d}),\tilde{\mathbf{m}}(\lambda,\mathbf{P},\mathbf{d})), \tag{4.10}$$

where $\mathbf{1}\in L(\mathbb{R}^2)$ denotes the identity, and $\mathbf{P}=(\mathbf{F},\mathbf{D})$, cf. (2.9). More explicitly, this reads

$$\begin{aligned} &div[W(\mathbf{F},\mathbf{D},\mathbf{d})\mathbf{1}-\mathbf{F}^TW_F(\mathbf{F},\mathbf{D},\mathbf{d})-\mathbf{D}^TW_D(\mathbf{F},\mathbf{D},\mathbf{d})] \\ &\qquad\qquad =\mathbf{F}^T\tilde{\mathbf{b}}(\lambda,\mathbf{P},\mathbf{d})+\mathbf{D}^T\tilde{\mathbf{m}}(\lambda,\mathbf{P},\mathbf{d})). \end{aligned}$$

Following an idea used in [13], we assume a stored-energy function of the form

$$W(\mathbf{P},\mathbf{d})=\Psi(\mathbf{P},\mathbf{d})+\Gamma(J), \tag{4.11}$$

where $\Psi:\Theta\to\mathbb{R}$ and $\Gamma:(0,\infty)\to\mathbb{R}$ are sufficiently smooth with

$$\Gamma\nearrow\infty \text{ as } J\searrow 0, \tag{4.12}$$

cf. (2.2). In what follows, we set $\tilde{\mathbf{F}}:=\mathbf{F}+\mathbf{d}\otimes\mathbf{e}_3$, and by (2.1), we have

$$J=\det\tilde{\mathbf{F}}=\det[\mathbf{F}+\mathbf{d}\otimes\mathbf{e}_3].$$

Specializing (4.10) to (4.11), we find

$$\begin{aligned} \nabla\Gamma(J)-div\{\Gamma'(J)\mathbf{F}^T[Cof\tilde{\mathbf{F}}]\mathbf{I}_o\}&=div\{\mathbf{F}^T\Psi_F(\cdot)+\mathbf{D}^T\Psi_D(\cdot)-\Psi(\cdot)\mathbf{1}\} \\ &\quad+\mathbf{F}^T\tilde{\mathbf{b}}(\lambda,\cdot)+\mathbf{D}^T\tilde{\mathbf{m}}(\lambda,\cdot), \end{aligned} \tag{4.13}$$

where the arguments $(\mathbf{F},\mathbf{D},\mathbf{d})\in\Theta$ are suppressed above.

**Lemma 4.1.** $\mathbf{F}^T[Cof\tilde{\mathbf{F}}]\mathbf{I}_o=(\det\tilde{\mathbf{F}})\mathbf{1}\Rightarrow div\{\Gamma'(J)\mathbf{F}^T[Cof\tilde{\mathbf{F}}]\mathbf{I}_o\}=\nabla\{J\Gamma'(J)\}$.

*Proof.* Observing that $\mathbf{F}=\tilde{\mathbf{F}}\mathbf{I}_o$, we find

$$\mathbf{F}^T[Cof\tilde{\mathbf{F}}]\mathbf{I}_o=\mathbf{I}_o^T\tilde{\mathbf{F}}^T[Cof\tilde{\mathbf{F}}^T]\mathbf{I}_o=\mathbf{I}_o^T\{(\det\tilde{\mathbf{F}})\mathbf{I}\}\mathbf{I}_o=(\det\tilde{\mathbf{F}})\mathbf{1}=J\mathbf{1},$$

where $\mathbf{I}\in L(\mathbb{R}^3)$ denotes the identity; the rest of the claim follows directly. □

Next, we define

$$\begin{aligned} &\Phi(J):=\Gamma(J)-J\Gamma'(J), \\ &\boldsymbol{\Xi}(\mathbf{P},\mathbf{d}):=\mathbf{P}^T\Psi_P(\mathbf{P},\mathbf{d})-\Psi(\mathbf{P},\mathbf{d})\mathbf{1}, \\ &\qquad\quad=\mathbf{F}^T\Psi_F(\mathbf{F},\mathbf{D},\mathbf{d})+\mathbf{D}^T\Psi_D(\mathbf{F},\mathbf{D},\mathbf{d})-\Psi(\mathbf{F},\mathbf{D},\mathbf{d})\mathbf{1} \end{aligned} \tag{4.14}$$

for all $(\mathbf{F},\mathbf{D},\mathbf{d})\in\Theta$, cf. (2.1). From Lemma 4.1 and (4.14), equation (4.13) then takes the convenient form

$$\nabla\Phi(J)=div\boldsymbol{\Xi}(\nabla\mathbf{f},\nabla\mathbf{d},\mathbf{d})+\nabla\mathbf{f}^T\tilde{\mathbf{b}}(\lambda,\nabla\mathbf{f},\nabla\mathbf{d},\mathbf{d})+\nabla\mathbf{d}^T\tilde{\mathbf{m}}(\lambda,\nabla\mathbf{f},\nabla\mathbf{d},\mathbf{d}). \tag{4.15}$$

**Theorem 4.2** *Assume that* (4.11) *satisfies the hypotheses of Theorem 4.1. In addition to* (4.12), *suppose that*

$$\Gamma'(d)<0 \text{ for } 0<d<d_o, \tag{4.16}$$

*for some* $d_o>0$, *and that there is a constant* $C_M>0$ *such that*

$$\left|\Xi_P(\mathbf{P},\mathbf{d})\,|+|\,\Xi_d(\mathbf{P},\mathbf{d})\right|\le C_M, \tag{4.17}$$

*for all* $(\mathbf{P},\mathbf{d})\in\Theta_M := \{(\mathbf{P},\mathbf{d})\in\Theta : |P|+|d|=|F|+|D|+|d|\le M\}$. *Then the solution branch* $\Sigma$ *is characterized by* (i) *and/or* (ii) *of Theorem 4.1.*

*Proof.* We use the same argument employed in [13]. Suppose that Theorem 4.1 (iii) is true, leading to (4.4) for a sequence $\{(\lambda_j,\mathbf{u}_j)\}\subset\Sigma$. We show that $\{\mathbf{u}_j\}$ is unbounded in $U$, arguing by contradiction: Suppose that $\{\mathbf{u}_j\}$ is bounded, and thus, (4.5) holds. From (4.14) and (4.15), we deduce that

$$\begin{aligned}\left\|\nabla\Phi(J_j)\right\|_\infty \le \left\|\Xi_P(\lambda_j,\mathbf{P}_j,\mathbf{d}_j)\{\nabla^2\mathbf{u}_j\}\right\|_\infty + \left\|\Xi_d(\lambda_j,\mathbf{P}_j,\mathbf{d}_j)[\nabla\mathbf{u}_j^2]\right\|_\infty \\ \left\|\mathbf{F}_j^T\tilde{\mathbf{b}}(\lambda_j,\mathbf{P}_j,\mathbf{d}_j)\right\|_\infty + \left\|\mathbf{D}_j^T\tilde{\mathbf{m}}(\lambda_j,\mathbf{P}_j,\mathbf{d}_j)\right\|_\infty,\end{aligned} \tag{4.18}$$

where $\mathbf{P}_j=(\mathbf{I}_o+\nabla\mathbf{u}_j^1,\nabla\mathbf{u}_j^2)$, $\mathbf{d}_j=\mathbf{e}_3+\mathbf{u}_j^2$, $J_j$ is as defined just after (4.4), and the quasi-linear notation $\{\nabla^2\mathbf{u}_j\}$, $[\nabla\mathbf{u}_j^2]$ are the same as those employed in (2.21), (2.22), respectively. In view of (2.23) and (4.17), each of the terms are the right side of (4.18) has a finite limit as $j\to\infty$. On the other hand, we claim that the limit of the term on the left side of the inequality grows without bound, which is a contradiction.

To see this, we use (4.6), (4.14)$_1$ and the presumed boundedness of $\{J_j\}$ to deduce that there are points $\mathbf{x}_o\in\bar{\Omega}\setminus\{\mathbf{x}_*\}$ such that $\Phi(J_j(\mathbf{x}_o))$ is bounded as $j\to\infty$. On the other hand, (4.12) and (4.16) imply that $\Phi(J_j(\mathbf{x}_*))\to\infty$ as $j\to\infty$, with $\mathbf{x}_*\in\bar{\Omega}$ as defined in (4.5). Finally, let $\hat{\mathbf{x}}(t)\subset\bar{\Omega}$ denote a smooth, unit-speed path, $0\le t\le 1$, with $\mathbf{x}_o=\hat{\mathbf{x}}(0)$ and $\mathbf{x}_*=\hat{\mathbf{x}}(1)$. The fundamental theorem of calculus and the mean value theorem yield

$$\Phi(J_j(\mathbf{x}_*))-\Phi(J_j(\mathbf{x}_0))\le \max_{t\in[0.1]}\left|\nabla_x\Phi(J_j(\hat{\mathbf{x}}(t))\right|.$$

Taking the limit as $j\to\infty$ above shows that the left side of (4.18) indeed grows without bound. Finally, since (4.10) can have no solution sequence $\{(\lambda_j,\mathbf{u}_j)\}$ with $\{\mathbf{u}_j\}\subset U$ bounded, we conclude that $\{\mathbf{u}_j\}$ and hence, $\Sigma$ must be unbounded. □

**Final Remarks**

Th results presented here are not special cases of those from our earlier works [10], [14] and [15]. The biggest difference comes from the role of strong ellipticity (2.11), which serves well for many of the required properties such as Lemma 3.2 and Propositions 4.1 and 4.2. However, material objectivity (2.3) brings the director $\mathbf{d}$ prominently into the picture, cf. (2.12)-(2.16). As already noted in Remark 2.1, the positivity assumption (2.1) is natural from the point of view of classical linear elasticity, and it plays a crucial role in Lemma 3.1. In nonlinear elasticity, the required properties associated with the deformation gradient are sufficient in the above-mentioned works. For instance, strong ellipticity alone is sufficient in [14] for global continuation in problems for homogeneous materials.

We are unable to eliminate alternative (ii) of Theorem 4.1, as done in [10] and [14]. By a uniqueness theorem from [18], a nontrivial solution in the presence of zero loading is not possible in those settings (cf. the discussion immediately after the statement of Theorem 4.1). On the other hand, our shell problem at the unloaded state is not comparable to the setup in [18]. If we imagine the single-director surface as a model for a slightly thickened elastic body in $\mathbb{R}^3$, then the lateral surfaces are not subjected to Dirichlet boundary conditions, as required in that work. As such, a uniqueness theorem of the type in [18] is not expected for the class of problems considered here.

As already pointed out in [10], the FPR degree [6] is convenient in global continuation problems. In particular, the spectral estimates in [14], [15] are not needed. However, the FPR degree requires a simply connected admissible space for global bifurcation problems. As in [14], [15], we must work in the set $\mathcal{O}$ as defined in (4.1), which is connected by definition - but not necessarily simply connected. As such, an oriented nonlinear Fredholm degree, such as that of [15], is apparently needed for bifurcation problems of nonlinear elasticity, cf. [12].

In the context of 3D nonlinear elasticity, the analogue of (4.15) was first uncovered in [3]; presuming additional regularity, energy minimizers are shown to satisfy the equation classically (in the absence of body forces). The same formulation was later employed in [14] to obtain a result like that of Theorem 4.2. We point out that Lemma 4.1 here is new. In any case, this corrects the omission of reference [3] in [14].

**Acknowledgements:** Thanks to P. Rosakis for a helpful conversation.

**Appendix**

*Proof of* (2.14): Differentiation of (3) yields

$$\begin{aligned} &W_F(\mathbf{QF},\mathbf{QD},\mathbf{Qd}) \equiv \mathbf{Q}W_F(\mathbf{F},\mathbf{D},\mathbf{d}),\\ &W_D(\mathbf{QF},\mathbf{QD},\mathbf{Qd}) \equiv \mathbf{Q}W_D(\mathbf{F},\mathbf{D},\mathbf{d}),\\ &W_d(\mathbf{QF},\mathbf{QD},\mathbf{Qd}) \equiv \mathbf{Q}W_d(\mathbf{F},\mathbf{D},\mathbf{d}),\ \ \forall \mathbf{Q} \in SO(3). \end{aligned}$$

Evaluating these at the reference stress-free reference configuration (12) then implies

$$\begin{aligned} &W_F(\mathbf{QI}_o,\mathbf{O},\mathbf{Qe}_3) \equiv \mathbf{O},\\ &W_D(\mathbf{QI}_o,\mathbf{O},\mathbf{Qe}_3) \equiv \mathbf{O},\\ &W_d(\mathbf{QI}_o,\mathbf{O},\mathbf{Qe}_3) \equiv \mathbf{0}. \end{aligned}$$

Substituting $\mathbf{Q} = \exp(\mathbf{W}t), \mathbf{W} \in Skew(\mathbb{R}^3), t \in \mathbb{R},$ above, differentiating with respect to $t,$ and evaluating at $t = 0,$ we find

$$W_{FF}(\mathbf{I}_o,\mathbf{O},\mathbf{e}_3)[\mathbf{WI}_o]+W_{Fd}(\mathbf{I}_o,\mathbf{O},\mathbf{e}_3)[\mathbf{We}_3]=\mathbf{O},$$
$$W_{DF}(\mathbf{I}_o,\mathbf{O},\mathbf{e}_3)[\mathbf{WI}_o]+W_{Dd}(\mathbf{I}_o,\mathbf{O},\mathbf{e}_3)[\mathbf{We}_3]=\mathbf{O},$$
$$W_{Fd}(\mathbf{I}_o,\mathbf{O},\mathbf{e}_3)[\mathbf{WI}_o]+W_{dd}(\mathbf{I}_o,\mathbf{O},\mathbf{e}_3)[\mathbf{We}_3]=\mathbf{0},$$

for all $\mathbf{W}\in Skew(\mathbb{R}^3)$. □

**References**

[1] Agmon, S., Douglis, A., Nirenberg, L.: Estimates near the boundary for solutions of elliptic partial differential equations satisfying general boundary conditions II, *Comm. Pure Appl. Math*., **17** (1964) 35–92.

[2] Antman, S.S., *Nonlinear Problems of Elasticity,* 2nd *Ed*., Springer-Verlag, 2005.

[3] Bauman, P., Phillips, D., Owen, N.C., Maximal smoothness of solutions to certain Euler-Lagrange equations from nonlinear elasticity, Proc. Royal Soc. Edinburgh 119A (1991), 241-263.

[4] Evans, L.C., *Partial Differential Equations, 2nd Ed.,* American Mathematical Society, Providence, RI, 2010.

[5] Fitzpatrick, P.M., Pejsachowicz, J., Orientation and Leray-Schauder theory for fully nonlinear elliptic boundary value problems, *Memoirs of the AMS,* Vol. 101-483, American Mathematical Society, Providence, RI, 1993.

[6] Fitzpatrick, P.M., Pejsachowicz, J., Rabier, P.J., The degree of proper $C^2$ Fredholm mappings I, *J. Reine Angew.Math*. **427** (1992) 1–33.

[7] Gilbarg, D., Trudinger, N.S., *Elliptic Partial Differential Equations of Second Order,* 2nd *Ed*., Springer-Verlag, 1983.

[8] Gurtin, M.E., *An Introduction to Continuum Mechanics,* Academic Press, 1981.

[9] Gurtin, M.E., Murdoch, A.I., A continuum theory of elastic material surfaces, *Arch. Rat. Mech. Anal.* **57** (1975) 291-323.

[10] Healey, T.J., Classical injective solutions in the large in incompressible nonlinear elasticity, *Arch. Rational Mech. Anal*. **232** (2018) 1207-1225.

[11] Healey, T.J., Kielhöfer, H., Symmetry and Nodal Properties in Global Bifurcation Analysis of Quasi-Linear Elliptic Equations, *Arch. Rat. Mech. Anal.* **113** (1991) 299-311.

[12] Healey, T.J., Montes-Pizarro, E.L.: Global bifurcation in nonlinear elasticity with an application to barreling states of cylindrical columns, *J. Elasticity,* **71,** 33–58, 2003.

[13] Healey, T.J., Nair, G.G., Energy minimizing configurations for single-director Cosserat shells, *J. Elasticity* **154** (2023) 569–578.

[14] Healey, T.J., Rosakis, P.: Unbounded branches of classical injective solutions in the forced displacement problem of nonlinear elastostatics, *J. Elasticity* **49** (1997) 65–78.

[15] Healey, T.J., Simpson, H.C., Global continuation in nonlinear elasticity, *Arch. Rational Mech. Anal*. **143** (1998) 1-28.

[16] Kato, T., *Perturbation Theory for Linear Operators*, 2nd Ed., Springer-Verlag, Berlin, 1980.

[17] Kielhöfer, H., *Bifurcation Theory*, 2nd *Ed*., Springer, New York, 2011.

[18] Knops, R.J., Stuart, C.A.: Quasiconvexity and uniqueness of equilibrium solutions in nonlinear elasticity, *Arch. Rational Mech. Anal*., **86**, 233-249, 1984.

[19] Li, Q., Healey, T.J., Stability boundaries for wrinkling in highly stretched elastic sheets, *J. Mech. Phys. Solids*, **97** (2016) 260-274.

[20] Nahgdi, P.M, Theory of Shells and Plates, *Handbuch der Physik,* Vol. VI/2, C. Truesdell, ed., Springer-Verlag Berlin (1972) 425-640.

[21] Nayyar, V., Ravi-Chandar, K., Huang, R., Stretch-induced stress patterns and wrinkles in hyperelastic thin sheets, *Int. J. Solids Struct*. **48** (2011) 3471–3483.

[22] Nugari, R., Further remarks on the Nemiskii operator in Hölder spaces, *Commentationes Mathematicae Universitatis Carolinae* **34** (1993) 89-95.

[23] Olver, P.J., *Applications of Lie Groups to Differential Equations,* Springer-Verlag, New York, 1986.

[24] Pejsachowicz, J., Rabier, P.J., A substitute for the Sard-Smale theorem in the $C^1$ case, *J. Anal. Math.* **76** (1998) 289-319.

[25] Rabinowitz, P.H.: Some global results for nonlinear eigenvalue problems, *J. Funct. Anal.*, **7** (1971) 487–513.

[26] Sokolnikoff, I.S., *Mathematical Theory of Elasticity, 2nd Ed.,* McGraw-Hill, 1956.

[27] Wloka, J., *Partial Differential Equations,* Cambridge University Press, Cambridge, 1987.